\documentclass{interact}

\usepackage{epstopdf}% To incorporate .eps illustrations using PDFLaTeX, etc.
\usepackage[caption=false]{subfig}% Support for small, `sub' figures and tables
\usepackage{dsfont}
\usepackage{tikz}
\usepackage{float}
\usepackage{booktabs}
\usepackage{verbatim}
\usepackage{stix}
\usepackage{algorithm}
\usepackage{algpseudocode}

\usepackage{amsmath}
\usepackage{tikz}
\usepackage{booktabs}
\usepackage{float}
\usepackage{geometry}

\usetikzlibrary{calc}

\usepackage[numbers,sort&compress]{natbib}% Citation support using natbib.sty
\bibpunct[, ]{[}{]}{,}{n}{,}{,}% Citation support using natbib.sty
\makeatletter% @ becomes a letter
\def\NAT@def@citea{\def\@citea{\NAT@separator}}% Suppress spaces between citations using natbib.sty
\makeatother% @ becomes a symbol again

\theoremstyle{plain}% Theorem-like structures provided by amsthm.sty
\newtheorem{theorem}{Theorem}[section]

\newtheorem{proposition}[theorem]{Proposition}

\theoremstyle{definition}

\theoremstyle{remark}

\begin{document}

\articletype{ARTICLE TEMPLATE}% Specify the article type or omit as appropriate

\title{Improving the efficiency of the ratio-of-uniforms method}

\author{
\name{Luis Rincón\textsuperscript{a}\thanks{CONTACT Luis Rincón. Email: lars@ciencias.unam.mx} and Emmanuel Delgadillo\textsuperscript{a}}
\affil{\textsuperscript{a}Departamento de Matemáticas, Facultad de Ciencias, UNAM}
}

\maketitle

\begin{abstract}
In 1977, Kinderman and Monahan showed that under certain conditions,
a random variable can be expressed as the ratio $V/U$,
where  $(U,V)$ is uniformly distributed over a bounded region  $S\subset{\mathds R}^2$.  Thus, by generating values of  $(U,V)$
and computing the indicated ratio, samples from the random variable can be obtained. This procedure is now known as the ratio-of-uniforms method
and relies heavily on our ability to produce uniform samples from $S$, for which
an acceptance-rejection method is customarily used.
In this paper, we propose three sampling schemes based on the polar representation of $S$. The first uses a covering consisting of triangles, the second uses a covering based on circular sectors,
while the third employs rectangular coverings of an associated angular density. Numerical examples illustrate that these schemes can achieve acceptance rates close to $100\%$.
\end{abstract}


\section{Introduction}
% ----------------------------

Efficient generation of values of a random variable is one of the central problems in stochastic simulation.
Monte Carlo methods
(see e.g.~\cite{Asmussen-Glynn-2007},
\cite{Ripley-1987},
\cite{Robert-Casella-2004}, and~\cite{Rubinstein-Kroese-2008})
consist of using those values to solve a variety of mathematical problems.
For generating values of a random variable,
one usually starts with one or two values from the uniform distribution on the interval $(0,1)$ and, through
some transformation, constructs a value from the target distribution.
For example, to obtain a random point uniformly distributed within a region $S\subset{\mathds R}^2$ of positive area and contained in a bounded rectangle $R=(a,b)\times(c,d)$, it is sufficient to generate a point inside the rectangle.
The point is accepted if it belongs to $S$ and rejected otherwise.
This can be done by independently choosing $u_1,u_2\in(0,1)$ and applying
the transformation $(a+(b-a)u_1,c+(d-c)u_2))$ gives a point in the rectangle $R$.
The efficiency (or acceptance rate) of such sampling is defined as the ratio $0<|S|/|R|\le1$,
where $|\cdot|$ denotes area.
It is clear that optimal efficiency is obtained when one chooses the rectangle (or some
other region $R$) of minimum area containing $S$ and for which a method to generate uniformly distributed points
is available.\\

The ratio-of-uniforms method is a general procedure to simulate a random variable
requiring only mild conditions on its distribution.
It was introduced in 1977 by
Kinderman and Monahan~\cite{Kinderman-Monahan-1977}.
In its univariate version, the method
establishes that if
$h(x)\ge0$ is a continuous, integrable, bounded function satisfying the conditions
\begin{equation}
\label{Conditions:for:method:quotient}
0<\displaystyle\int_{-\infty}^{\infty}h(x)\,dx<\infty\quad\mbox{and}\quad
\displaystyle\sup_{x}\,|x|\sqrt{h(x)}<\infty,
\end{equation}
then the region
\begin{equation}
\label{Region:S:cartesian:coord}
S=\{(u,v)\in{\mathds R}^2:\ 0<u<\sqrt{h(v/u)},\ v/u\in\mathrm{Supp}(h)\}\subset{\mathds R}^{2},
\end{equation}
is bounded, has area $|S|=(1/2)\,\int_{-\infty}^{\infty}h(x)\,dx$,
and has the property that if  $(u,v)$ is a random point uniformly distributed on $S$,
then $v/u$ is a sample from the density
$f(x)=h(x)/\int_{-\infty}^{\infty}h(x)dx$.\\

A proof of this seemingly enigmatic result can be found in the original article~\cite{Kinderman-Monahan-1977}.
It can also be found  in~\cite{Dagpunar-2007},
\cite{Gentle-2003},
\cite{Nelson-2013},
\cite{Press-etal},
\cite{Rubinstein-Kroese-2008} and in most modern books on stochastic simulation.
To generate points in $S$ using the acceptance-rejection method, one usually encloses $S$ within a rectangle
of the smallest possible area in order to maximize the acceptance rate. The use of
rectangles makes the uniform sampling straightforward with the aid of a computer.
From~\eqref{Region:S:cartesian:coord}, it can be shown that the cartesian coordinates $(u,v)$ 
of the outer boundary
of the region $S$, wherever $h(x)>0$, can be parametrized by the values $x$ in the support of $h$ by
$(u,v)=(\sqrt{h(x)},x\,\sqrt{h(x)})$.
This makes it possible to determine the rectangle of minimum area (with sides parallel to the coordinate axes)
containing $S$ and, in this way,
increase the efficiency of the sampling procedure. The resulting rectangle is
\begin{equation}
\label{Rectángulo-óptimo-método-estándar-unidim}
R={\big(}0,\sup_x\sqrt{h(x)}\,{\big)}\times{\big(}\inf_x x\,\sqrt{h(x)},\ \sup_x x\,\sqrt{h(x)}\,{\big)}.
\end{equation}

In Bayesian analysis, prior distributions are often specified
without the multiplicative normalization constant because these constants may be difficult to calculate.
The ratio-of-uniforms method is particularly useful in simulating values from such
distributions as the method does not require $h(x)$ to be a density.\\

Since its introduction, the method has been the subject of several improvements and extensions.
For example, Barabesi (1993)~\cite{Barabesi-1993} provides a broad review of the method,
Jones and Lunn (1996)~\cite{Jones-Lunn-1996} present a pedagogical exposition,
emphasizing (quite appropriately) that transformations are a fundamental tool
in the various methods for generating random variables. Indeed, the ratio-of-uniforms method uses the
one-to-one transformation $\Phi(x,y)=(\sqrt{y},x\sqrt{y})$ to take $S_0=\{(x,y):x\in\mbox{Supp}(h), 0<y<h(x)\}$ into $S$
such that if $(U,V)\sim\mbox{Unif}(S)$ then $\Phi^{-1}(U,V)\sim\mbox{Unif}(S_0)$.
In Li et al. (2023)~\cite{li2023objectivebayesiananalysisgeneralized}, the authors carry out Bayesian estimation
for the two parameters of a generalized exponential distribution. To do so, they employ the
$r$-parameter ratio-of-uniforms method to obtain independent posterior samples from certain marginal distributions.
Perez et al. (2008)~\cite{perez2008efficient} proposed using hyperellipsoids
to enclose the region $S$ and thereby achieve
greater efficiency in uniform sampling over $S$, particularly in the multivariate setting.
On the other hand, Leydold (2001)~\cite{Leydold2001} and Leydold (2003)~\cite{Leydold2003} discuss several proposals
for the simulation of discrete distributions based on the ratio-of-uniforms method.
See Press et al.~\cite{Press-etal} for some worked-out examples in the case of discrete random variables.\\

The aim of this paper is to propose three simple sampling schemes that improve the acceptance rate of the ratio-of-uniforms
method for univariate distributions. All three schemes are based on the polar representation of the associated region $S$.
The first employs a triangular covering, the second uses circular sectors, and the third is based on a rectangular covering
of an associated angular density. We establish convergence of the acceptance rates to one and show that the second and
third schemes have identical acceptance rates for a given partition.\\

% -----------------------------------------
\section{Polar representation}
% -----------------------------------------
Using polar coordinates, the region $S$ given in~\eqref{Region:S:cartesian:coord} can be written as
\begin{equation}
\label{Region:S:polar:coord}
S=\{(r,\theta):\ 0<r<\rho(\theta),\ -\pi/2<\theta<\pi/2,\},
\end{equation}
where the radial function is
\begin{equation}
\label{Radial:function}
\rho(\theta)=\frac{\sqrt{h(\tan\theta)}}{\cos\theta},\quad -\pi/2<\theta<\pi/2.
\end{equation}
Assume $\rho(\theta)$ admits a continuous extension to the closed interval
$[-\pi/2,\pi/2]$. Then $\rho(\theta)$ is bounded on $[-\pi/2,\pi/2]$ and a change of variable yields
\begin{displaymath}
\int_{-\pi/2}^{\pi/2}\rho^2(\theta)\,d\theta=\int_{-\infty}^{\infty} h(x)\,dx.
\end{displaymath}

% -----------------------------------------
\section{A triangular covering}
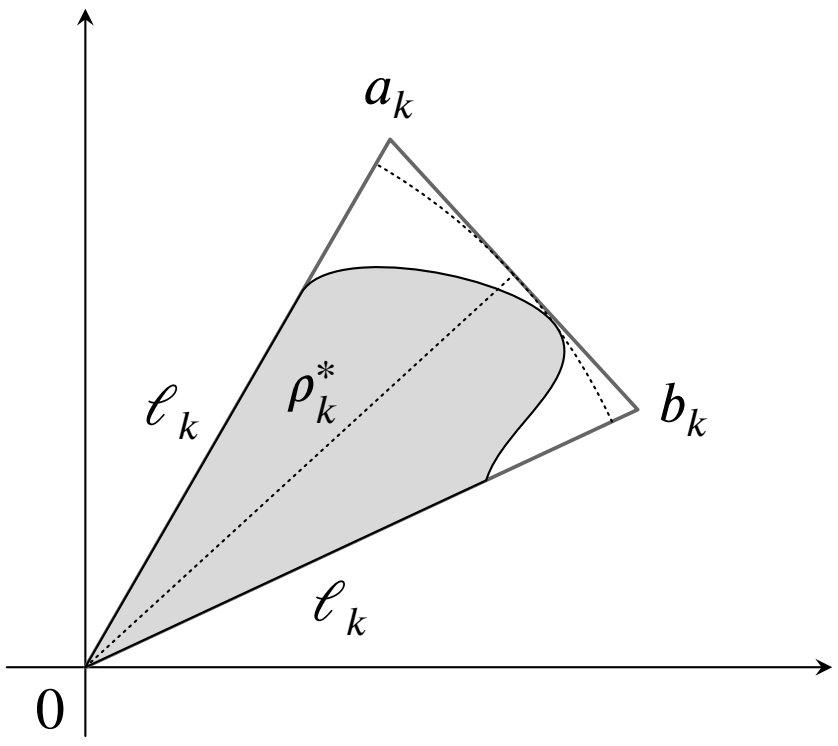
\begin{figure}[H]
\centering

\begin{tikzpicture}[
    scale=1.3,
    >=stealth,
    line cap=round,
    line join=round
]

% Parameters
\pgfmathsetmacro{\L}{3.10}
\pgfmathsetmacro{\angleA}{60}
\pgfmathsetmacro{\angleB}{25}
\pgfmathsetmacro{\angleM}{(\angleA+\angleB)/2}

% Angle where the Bezier curve reaches its maximum
%\pgfmathsetmacro{\angleT}{49}
\pgfmathsetmacro{\angleT}{36}

% Radius of the circle tangent to AB
\pgfmathsetmacro{\R}{
    \L*cos((\angleA-\angleB)/2)
}

% Main points
\coordinate (O) at (0,0);
\coordinate (A) at (\angleA:\L);
\coordinate (B) at (\angleB:\L);
\coordinate (M) at (\angleM:\R);

% Point of contact with the circular arc
\coordinate (T) at (\angleT:\R);

% Coordinate axes
\draw[->] (-0.4,0) -- (3.8,0);
\draw[->] (0,-0.35) -- (0,3.35);

% Triangle
\draw[line width=0.7pt, draw=gray!120] (O) -- (A) -- (B) -- cycle;

% Dotted circular arc
\draw[densely dotted]
    (\angleB:\R)
    arc[start angle=\angleB,
        end angle=\angleA,
        radius=\R];

% Bezier curve with maximum radius at T
% Tangent direction at T is angleT - 90 degrees

\coordinate (P) at (\angleA:2.20);
\coordinate (Q) at (\angleB:2.25);

% Region bounded by the exterior radii
% and the Bezier curve

\filldraw[fill=gray!30,draw=black]
    (O) -- (P)
    .. controls (60:2.50)
             and ($(T)+(\angleT+90:0.32)$) ..
    (T)
    .. controls ($(T)+(\angleT-90:0.32)$)
             and (30:2.45) ..
    (Q)
    -- cycle;
    
    % Single dotted radius: angular bisector
\draw[densely dotted]
    (O) -- (M);

% Labels
\node[below left] at (O) {$0$};
\node[above] at (A) {$a_k$};
\node[right] at (B) {$b_k$};

\node[left] at (1.4,1.4) {$\rho_k^*$};
\node[left] at (0.7,1.3) {$\ell_k$};
\node[below] at (1.3,0.55) {$\ell_k$};

\end{tikzpicture}

   \caption{Construction of isosceles triangle $\Delta_k$ with height $\rho_k^*$.}
    \label{Triangle-construction}
\end{figure}

% -----------------------------------------
The polar representation is here used to construct a covering for $S$ using a collection of triangles.
Let $n\ge2$ be an integer and consider a uniform partition $\{\theta_0<\theta_1<\cdots<\theta_n\}$ of
$(-\pi/2,\pi/2)$ given by $\theta_k=-\pi/2+k\pi/n$, $k=0,1,\ldots,n$. For each $k=1,\ldots,n$, let $\rho_k^*$ be the
maximum value of the radial function $\rho(\theta)$ on $[\theta_{k-1},\theta_k]$. Using  $\rho_k^*$ as a constant radius,
a circular sector can be drawn from angle $\theta_{k-1}$ to angle $\theta_k$. This is shown in Figure~\ref{Triangle-construction}.
An isosceles triangle $\Delta_k$ can be defined  with one vertex at the origin (apex) and the other two vertices at
\begin{align*}
a_k&=\ell_k(\cos\theta_{k-1},\sin\theta_{k-1}),\\
b_k&=\ell_k(\cos\theta_{k},\sin\theta_{k}),
\end{align*}
where the two equal sides have length 
\begin{displaymath}
\ell_k=\frac{\rho_k^*}{\cos(\pi/(2n))}.
\end{displaymath}

The height $\rho_k^*$ bisects the vertex angle \(\theta_k-\theta_{k-1}=\pi/n$.
Thus,  a covering for $S$ can be defined as $\Delta:=\bigcup_{k=1}^n\Delta_k$,
where it can be shown that $\Delta_k$ has area $|\Delta_k|=(\rho_k^*)^2\tan(\pi/(2n))$. The area of the
covering is
\begin{displaymath}
|\Delta|=\tan\big(\frac{\pi}{2n}\big)\sum_{k=1}^n(\rho_k^*)^2.
\end{displaymath}

A uniformly distributed random point in \(\Delta_k\) can be selected using the formula
$ua+vb$, where $a$ and $b$ are the two nonzero vertices and $u$ and $v$
are two independent $\mbox{Unif}(0,1)$ values. Here $u$ is replaced by $1-u$, and $v$ is
replaced by $1-v$ when the chosen pair $(u,v)$ satisfies $u+v>1$.
To obtain a uniform point from $\Delta$, first a triangle $\Delta_k$ is selected with
probability $|\Delta_k|/|\Delta|$. Having chosen one $\Delta_k$,
a point from it is selected as explained above.
The result is a uniformly distributed point $(u,v)\in\Delta$, which is accepted when it belongs to $S$, and is rejected otherwise.
When accepted, $v/u$ is a sample from the density $h(x)/\int_{-\infty}^{\infty}h(x)dx$. The rate of acceptance is
\begin{equation}
\label{Acceptance-rate-triangular}
A_n(\Delta):=\frac{|S|}{|\Delta|}=\frac{\int_{-\infty}^{\infty}h(x)\,dx}{2\tan\big({\pi}/{(2n)}\big)\sum_{k=1}^n(\rho_k^*)^2}.
\end{equation}

\begin{proposition}
The rate of acceptance $|S|/|\Delta|$ using a triangular covering converges to $1$ as the partition size $n$ tends to infinity. 
\end{proposition}

\begin{proof}
Since $h(x)$ is bounded and continuous, so is $\rho^2(\theta)$.
By Darboux's theorem~\cite{bartlesherbert2011} on integration, the weighted average of the maxima
(upper Darboux sum)
converges to the Riemann integral of $\rho^2(\theta)$, i.e.
\begin{align*}
\lim_{n\to\infty}2\tan\big(\frac{\pi}{2n}\big)\sum_{k=1}^n(\rho_k^*)^2
&=\lim_{n\to\infty} \frac{2n}{\pi}\tan\big(\frac{\pi}{2n}\big)
\cdot\lim_{n\to\infty}\frac{\pi}{n}\sum_{k=1}^n(\rho_k^*)^2\\
&=\lim_{n\to\infty}\frac{\pi}{n}\sum_{k=1}^n(\rho_k^*)^2\\
&=\int_{-\pi/2}^{\pi/2}\rho^2(\theta)\,d\theta\\
&=\int_{-\infty}^{\infty}h(x)\,dx.
\end{align*}
\end{proof}

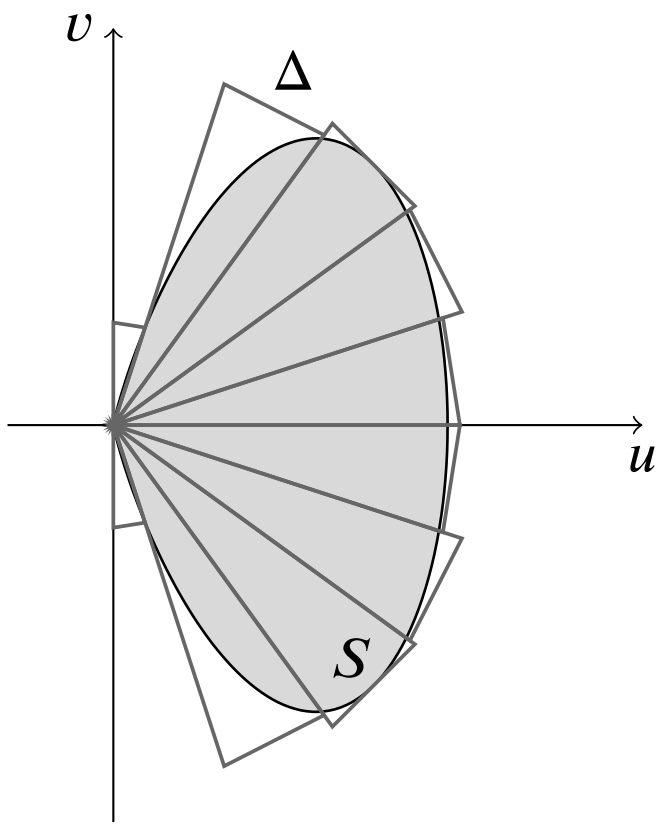
\begin{figure}[H]
    \centering
    \begin{tikzpicture}[xscale=3.5, yscale=3.5]

        % Ratio-of-uniforms region for the standard normal density
        \path[
            draw=black,
            line width=0.5pt,
            fill=gray!30
        ]
            (0,0)
            -- plot[
                domain=-8:8,
                samples=401,
                variable=\t
            ]
            ({exp(-\t*\t/4)/(2*pi)^0.25},
             {\t*exp(-\t*\t/4)/(2*pi)^0.25})
            -- cycle;

        % Cartesian axes
        \draw[->] (0,-0.75) -- (0,0.75)
            node[left] {$v$};
        \draw[->] (-0.2,0) -- (1,0)
            node[below] {$u$};
            
             \node%[fill=gray!30,inner sep=1pt]
            at (0.34,0.67) {$\Delta$};
            
        \node[fill=gray!30,inner sep=1pt]
            at (0.45,-0.44) {$S$};

        % Triangular sections
        \begin{scope}[line width=0.7pt, draw=gray!120]
            \draw (0,0) -- (0,-0.1938)
                -- (0.0599,-0.1843) -- cycle;

            \draw (0,0) -- (0.2094,-0.6444)
                -- (0.3982,-0.5481) -- cycle;

            \draw (0,0) -- (0.414,-0.5698)
                -- (0.5698,-0.414) -- cycle;

            \draw (0,0) -- (0.5604,-0.4072)
                -- (0.6588,-0.2141) -- cycle;

            \draw (0,0) -- (0.6228,-0.2024)
                -- (0.6549,0) -- cycle;

            \draw (0,0) -- (0.6549,0)
                -- (0.6228,0.2024) -- cycle;

            \draw (0,0) -- (0.6588,0.2141)
                -- (0.5604,0.4072) -- cycle;

            \draw (0,0) -- (0.5698,0.414)
                -- (0.414,0.5698) -- cycle;

            \draw (0,0) -- (0.3982,0.5481)
                -- (0.2094,0.6444) -- cycle;

            \draw (0,0) -- (0.0599,0.1843)
                -- (0,0.1938) -- cycle;
        \end{scope}

        % Table to the right of the plot
        \node[anchor=west,inner sep=0pt] at (1.5,0) {%
            \renewcommand{\arraystretch}{1}%
            \begin{tabular}{cc}
                \toprule
                \shortstack{Partition\\
                            size $n$}
                &
                \shortstack{Acceptance rate\\
                            $|S|/|\Delta|\ (\%)$}
                \\
                \midrule
                10    & 85.09\% \\
                50   & 96.52\% \\
                100   & 98.24\% \\
                500  & 99.64\% \\
                1000  & 99.82\% \\
                1500  & 99.88\% \\
                10000 & 99.99\% \\
                \bottomrule
            \end{tabular}%
        };

    \end{tikzpicture}
    \caption{A triangular covering ($n=10$) and acceptance rates of sampling
             for the standard normal distribution.}
    \label{Normal:triangular:partition}
\end{figure}

As an example, consider the kernel $h(x)=\exp\{-x^2/2\}$
of the standard normal density. The region $S$ is given by
\begin{displaymath}
S=\{(u,v):0<u<1,\quad v^2<-4u^2\ln u\},
\end{displaymath}
with radial function $\rho(\theta)=(1/\cos\theta)\,\exp\{-(\tan^2\theta)/4\}$.
The plot on the left of Figure~\ref{Normal:triangular:partition} shows a triangular covering $\Delta$ for $S$
and the table on the right shows how the acceptance rate increases as
the partition size $n$ grows.

% ----------------------------------------------------------
\section{A circular-sector covering}
% ----------------------------------------------------------
Using the square root method~\cite{Press-etal}, a uniform sample from the circular sector shown in Figure~\ref{Triangle-construction}
can be obtained
by first choosing an angle uniformly between $\theta_{k-1}$ and $\theta_k$, and then choosing a radius of length
$\rho_k^*\sqrt{u}$, where $u$ is a sample from $\mbox{Unif}(0,1)$. Denote this circular sector by $C_k$.
%which is shown in Figure~\ref{Triangle-construction}.
Then  a smaller covering for $S$ can be defined as $C:=\bigcup_{k=1}^nC_k$,
where $|C_k|=(1/2)(\rho_k^*)^2(\pi/n)$. The area of the
covering is
\begin{displaymath}
|C|=\big(\frac{\pi}{2n}\big)\sum_{k=1}^n(\rho_k^*)^2.
\end{displaymath}
Thus, the  acceptance rate  is
\begin{equation}
\label{Acceptance-rate-circular}
A_n(C):=\frac{|S|}{|C|}=\frac{\int_{-\infty}^{\infty}h(x)\,dx}{({\pi}/{n})\sum_{k=1}^n(\rho_k^*)^2}.
\end{equation}
Since $x<\tan x$, for $0<x<\pi/2$, a comparison between the acceptance rates~\eqref{Acceptance-rate-triangular} and~\eqref{Acceptance-rate-circular}
gives that the circular-sector covering has a higher rate,
\begin{align*}
\frac{|S|}{|\Delta|}&=\frac{\int_{-\infty}^{\infty}h(x)\,dx}{2\tan\big({\pi}/{(2n)}\big)\sum_{k=1}^n(\rho_k^*)^2}\\
&=\frac{\pi/(2n)}{\tan\big({\pi}/{(2n)}\big)}\cdot\frac{\int_{-\infty}^{\infty}h(x)\,dx}{({\pi}/{n})\sum_{k=1}^n(\rho_k^*)^2}\\
&=\frac{\pi/(2n)}{\tan\big({\pi}/{(2n)}\big)}\cdot\frac{|S|}{|C|}\\
&\le\frac{|S|}{|C|}.
\end{align*}

The plot on the left of Figure~\ref{Normal:circular:partition} shows a circular-sector covering $C$ for $S$
and the table on the right shows how the acceptance rate increases as the partition size $n$ grows.
As expected, the acceptance rates are slightly higher than those for the triangular covering.
Values are displayed up to two decimal places.

% Radial function for the standard normal
% Angles are measured in degrees.
%
% This function is evaluated only for
% angles strictly between -90 and 90.

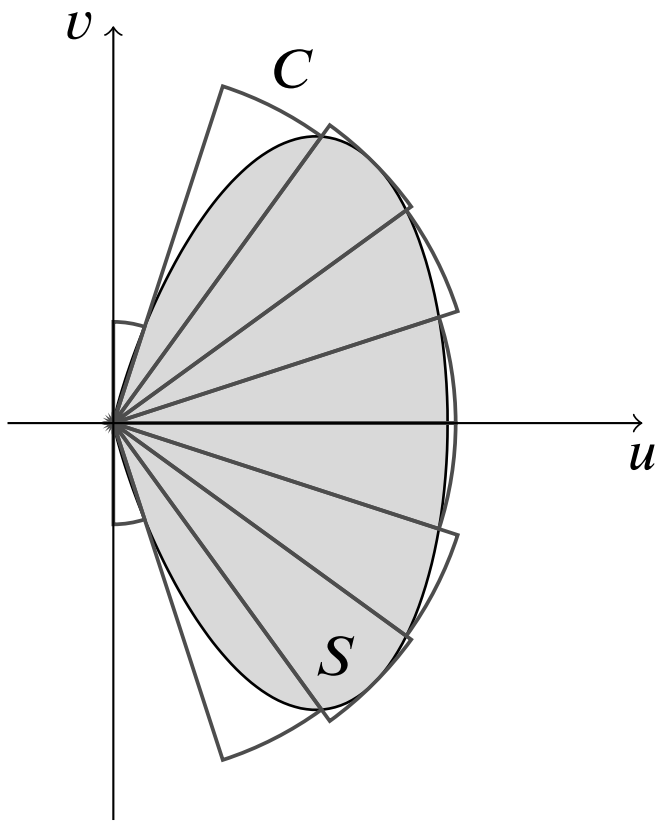
\pgfmathdeclarefunction{radial}{1}{%
    \pgfmathparse{
        sqrt(1+tan(#1)^2)*
        exp(-tan(#1)^2/4)/
        (2*pi)^0.25
    }%
}

\begin{figure}[H]
    \centering

    \begin{tikzpicture}[xscale=3.5,yscale=3.5]

        %------------------------------------
        % Region S
        %------------------------------------

        \path[
            draw=black,
            line width=0.5pt,
            fill=gray!30
        ]
            (0,0)
            -- plot[
                domain=-8:8,
                samples=401,
                variable=\t
            ]
            ({exp(-\t*\t/4)/(2*pi)^0.25},
             {\t*exp(-\t*\t/4)/(2*pi)^0.25})
            -- cycle;

%------------------------------------
% Circular sector covering
%------------------------------------

\pgfmathtruncatemacro{\n}{10}
\pgfmathsetmacro{\dtheta}{180/\n}

\begin{scope}[
    line width=0.7pt,
    draw=black!70
]

\foreach \k in {0,...,9} {

    % Angular interval
    \pgfmathsetmacro{\thetaA}{
        -90+\k*\dtheta
    }
    \pgfmathsetmacro{\thetaB}{
        \thetaA+\dtheta
    }

    % Endpoint radii
    \pgfmathsetmacro{\Ra}{0}
    \pgfmathsetmacro{\Rb}{0}

    \ifdim \thetaA pt>-89.999pt
        \pgfmathsetmacro{\Ra}{
            radial(\thetaA)
        }
    \fi

    \ifdim \thetaB pt<89.999pt
        \pgfmathsetmacro{\Rb}{
            radial(\thetaB)
        }
    \fi

    % Initial maximum
    \pgfmathsetmacro{\Rmax}{
        max(\Ra,\Rb)
    }

    % Critical point -45 degrees
    \ifdim \thetaA pt<-45pt
        \ifdim \thetaB pt>-45pt
            \pgfmathsetmacro{\Rmax}{
                max(\Rmax,radial(-45))
            }
        \fi
    \fi

    % Critical point 0 degrees
    \ifdim \thetaA pt<0pt
        \ifdim \thetaB pt>0pt
            \pgfmathsetmacro{\Rmax}{
                max(\Rmax,radial(0))
            }
        \fi
    \fi

    % Critical point 45 degrees
    \ifdim \thetaA pt<45pt
        \ifdim \thetaB pt>45pt
            \pgfmathsetmacro{\Rmax}{
                max(\Rmax,radial(45))
            }
        \fi
    \fi

    % Draw the covering sector
    \draw
        (0,0)
        -- (\thetaA:\Rmax)
        arc[
            start angle=\thetaA,
            end angle=\thetaB,
            radius=\Rmax
        ]
        -- cycle;
}

\end{scope}

        %------------------------------------
        % Coordinate axes
        %------------------------------------

        \draw[->]
            (0,-0.75) -- (0,0.75)
            node[left] {$v$};

        \draw[->]
            (-0.2,0) -- (1,0)
            node[below] {$u$};

        %------------------------------------
        % Labels
        %------------------------------------

        \node at (0.34,0.67) {$C$};

        \node[
            fill=gray!30,
            inner sep=1pt
        ] at (0.42,-0.44) {$S$};

        %------------------------------------
        % Table to the right
        %------------------------------------

        \node[
            anchor=west,
            inner sep=0pt
        ] at (1.5,0) {%

            \renewcommand{\arraystretch}{1}%

            \begin{tabular}{cc}
                \toprule

                \shortstack{
                    Partition\\
                    size $n$
                }
                &
                \shortstack{
                    Acceptance rate\\
                    $|S|/|C|\ (\%)$
                }
                \\

                \midrule

                10    & 85.80\% \\
                50   & 96.55\% \\
                100   & 98.25\% \\
                500  & 99.64\% \\
                1000  & 99.82\% \\
                1500  & 99.88\% \\
                10000 & 99.99\% \\

                \bottomrule
            \end{tabular}%
        };

    \end{tikzpicture}

    \caption{
        A circular-sector covering ($n=10$)
        and acceptance rates of sampling
        for the standard normal distribution.
    }

    \label{Normal:circular:partition}

\end{figure}

\begin{proposition}
The rate of acceptance $|S|/|C|$ using a circular-sector covering converges to $1$ as the partition size $n$ tends to infinity. 
\end{proposition}

The proof follows the same argument as in the triangular covering.
By Darboux's theorem, the weighted average of the maxima
converges to the Riemann integral of $\rho^2(\theta)$, i.e.
\begin{displaymath}
\lim_{n\to\infty}\frac{\pi}{n}\sum_{k=1}^n(\rho_k^*)^2
=\int_{-\pi/2}^{\pi/2}\rho^2(\theta)\,d\theta
=\int_{-\infty}^{\infty}h(x)\,dx.
\end{displaymath}

% ---------------------------------------------------------------------------------
\section{A rectangular covering for the angular density}
% ---------------------------------------------------------------------------------

This is another simple simulation technique based on the ratio-of-uniforms method and the following observation.

\begin{proposition}
\label{Angular:density}
Let $h(x)$ be a function satisfying the conditions of the ratio-of-uniforms method
with normalization constant $c=1/\int_{-\infty}^{\infty}h(x)dx$.
Let $\rho(\theta)={\sqrt{h(\tan\theta)}}/{\cos\theta}$ be the radial function.
Then
\begin{enumerate}
\item If $X$ has density $ch(x)$ then $Y=\arctan X$ has density
$c\rho^2(\theta)$.
\item If $Y$ has density $c\rho^2(\theta)$ then $X=\tan Y$ has density $ch(x)$.
\end{enumerate}
\end{proposition}

\begin{proof}
Both results follow from the change-of-variables formula.
For $(1)$,
\begin{displaymath}
f_Y(\theta)=\frac{ch(\tan \theta)}{\cos^2 \theta}
=c\rho^2(\theta).
\end{displaymath}
For $(2)$,
\begin{displaymath}
f_X(x)=\frac{c\rho^2(\arctan x)}{1+x^2}
=\frac{ch(x)}{(\cos(\arctan x))^2}\,\frac{1}{1+x^2}
=ch(x).
\end{displaymath}
\end{proof}

Thus, since $c\rho^2(\theta)$ is a bounded density defined
over a bounded domain, the acceptance-rejection method can be applied to obtain values of $Y$, from which values of $X$ can be derived. We will refer to $c\rho^2(\theta)$
as the angular density associated with $h(x)$.
In this scheme the region $S_a$ to sample from
is defined as
\begin{equation}
\label{Angular:set:S}
S_a=\{(\theta,v):-\pi/2<\theta<\pi/2,\quad 0<v<c\rho^2(\theta)\},
\end{equation}
with area $|S_a|=1$.
%\int_{-\infty}^{\infty}h(x)dx$.
As before, let $\{\theta_0<\theta_1<\cdots<\theta_n\}$ be a finite partition of
$(-\pi/2,\pi/2)$ given by $\theta_k=-\pi/2+k\pi/n$, $k=0,1,\ldots,n$. 
For each $k=1,\ldots,n$, let $m_k^2$ be the
maximum value of $c\rho^2(\theta)$ on $[\theta_{k-1},\theta_k]$
and define the rectangle
\begin{equation}
H_k=\{(\theta,v):\theta_{k-1}<\theta<\theta_k,\ \ 0<v<m_k^2\}.
%\vrectangle
\end{equation}
as in an upper Riemann sum, where $|H_k|=(\pi/n)m_k^2$.
Thus,  a covering for $S_a$ can be defined as $H:=\bigcup_{k=1}^nH_k$, with area $|H|=(\pi/n)\sum_{k=1}^nm_k^2$.
To generate a uniform point from $S_a$, first a rectangle $H_k$ is selected.
Again, it is important  here that the rectangle
$H_k$ be chosen according to its weight $|H_k|/|H|$. Having chosen one $H_k$,
a point from it is then selected in the usual way.
The result is a point $(\theta,v)$ which is accepted when it belongs to $S_a$ and
is rejected otherwise. When accepted, the first coordinate $\theta$ is a sample from the angular density $f_Y(\theta)=c\rho^2(\theta)$
and hence, $\tan\theta$ is a sample from the density $h(x)/\int_{-\infty}^{\infty}h(x)dx$. The acceptance rate is
\begin{displaymath}
A_n(H):=\frac{|S_a|}{|H|}=\frac{1}{({\pi}/{n})\sum_{k=1}^nm_k^2}.
\end{displaymath}

\begin{proposition}
The rate of acceptance $|S_a|/|H|$ of the rectangular covering of the angular density $c\rho^2(\theta)$
converges to $1$ as the partition size $n$ tends to infinity. 
\end{proposition}

As before, the proof follows the same argument as in the two previous coverings.
By Darboux's theorem, the weighted average of the maxima
converges to the Riemann integral of $c\rho^2(\theta)$, i.e.
\begin{displaymath}
\lim_{n\to\infty}\frac{\pi}{n}\sum_{k=1}^nm_k^2
=\int_{-\pi/2}^{\pi/2}c\rho^2(\theta)\,d\theta
=1.
\end{displaymath}

The following result shows that the rectangular covering of the angular density is equivalent to the
circular-sector covering of $S$ in terms of efficiency.

\begin{proposition}
Let $h(x)$ be a function satisfying the conditions of the ratio-of-uniforms method with
normalization constant $c=1/\int_{-\infty}^{\infty}h(x)dx$, with region $S$ and radial function $\rho(\theta)$.
The rate of acceptance of the rectangular covering of the angular density $c\rho^2(\theta)$
coincides with the rate of acceptance of the circular-sector covering of $S$.
\end{proposition}

\begin{proof}
\begin{displaymath}
\frac{|S|}{|C|}=\frac{\int_{-\infty}^{\infty}h(x)\,dx}{({\pi}/{n})\sum_{k=1}^n(\rho_k^*)^2}
=\frac{1}{({\pi}/{n})\sum_{k=1}^nc(\rho_k^*)^2}
=\frac{1}{({\pi}/{n})\sum_{k=1}^nm_k^2}
=\frac{|S_a|}{|H|}.
\end{displaymath}
\end{proof}

As an example, consider again the kernel $h(x)=\exp\{-x^2/2\}$ of the
standard normal density, $c=1/\sqrt{2\pi}$. Assume $X$ has $\mbox{N}(0,1)$ density. From Proposition~\ref{Angular:density}, $Y=\arctan X$ has density
\begin{displaymath}
f_Y(\theta)=\frac{1}{\sqrt{2\pi}}\,\frac{1}{\cos^2\theta}\,\exp\{-(\tan\theta)^2/2\},
\quad -\pi/2<\theta<\pi/2.
\end{displaymath}

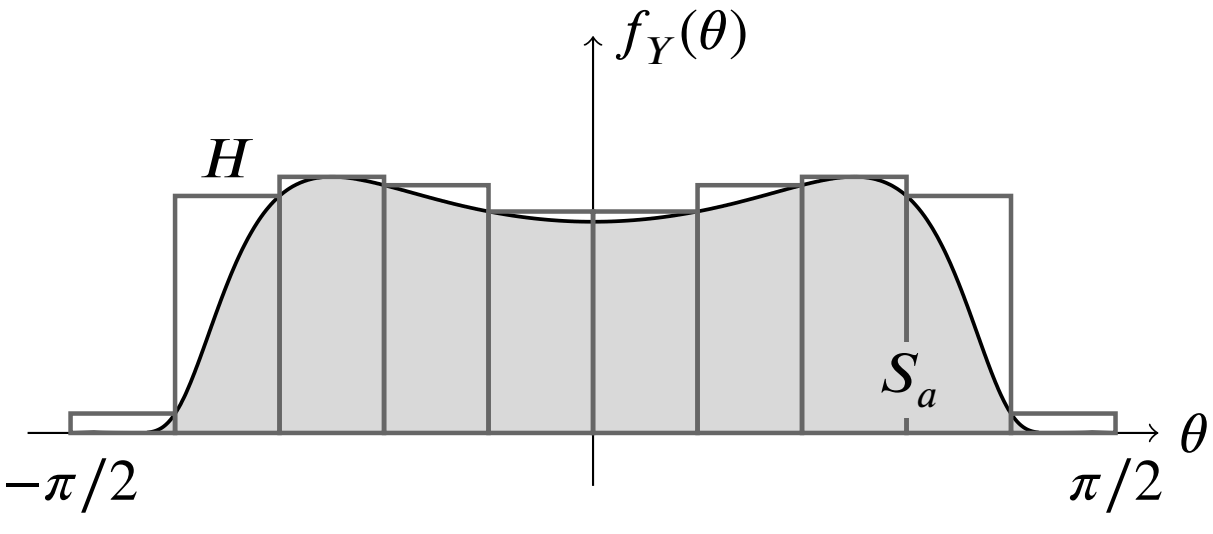
\begin{figure}[H]
    \centering
    \begin{tikzpicture}[xscale=2.2,yscale=3.5]
    
            % Cartesian axes
        \draw[->] (0,-0.1) -- (0,0.75)
            node[right] {$f_Y(\theta)$};
        \draw[->] (-1.7,0) -- (1.7,0)
            node[right] {$\theta$};

        % Graph of rho^2
\path[
    draw=black,
    line width=0.7pt,
    fill=gray!30,
]
    (-pi/2,0)
    -- plot[
        domain=-1.5:1.5,
        samples=1001,
        variable=\t
    ]
    (\t,{
        exp(-0.5*tan(deg(\t))^2)
        /(sqrt(2*pi)*cos(deg(\t))^2)
    })
    -- (pi/2,0)
    -- cycle;

        % Rectangle coverings and their labels
        \begin{scope}[line width=0.85pt, draw=gray!120]
            \draw (-1.5708,0) rectangle (-1.2566,0.0367);
            %\node at (-1.415,0.10) {\small $1$};

            \draw (-1.2566,0) rectangle (-0.9425,0.4478);
            %\node at (-1.1,0.2) {\small $2$};

            \draw (-0.9425,0) rectangle (-0.6283,0.4839);
            %\node at (-0.8,0.25) {\small $3$};

            \draw (-0.6283,0) rectangle (-0.3142,0.4681);
            %\node at (-0.48,0.23) {\small $4$};

            \draw (-0.3142,0) rectangle (0,0.4184);
            %\node at (-0.16,0.18) {\small $5$};

            \draw (0,0) rectangle (0.3142,0.4184);
            %\node at (0.16,0.18) {\small $6$};

            \draw (0.3142,0) rectangle (0.6283,0.4681);
            %\node at (0.48,0.23) {\small $7$};

            \draw (0.6283,0) rectangle (0.9425,0.4839);
            %\node at (0.8,0.25) {\small $8$};

            \draw (0.9425,0) rectangle (1.2566,0.4478);
            %\node at (1.1,0.2) {\small $9$};

            \draw (1.2566,0) rectangle (1.5708,0.0367);
            %\node at (1.415,0.10) {\small $10$};
        \end{scope}

            % Horizontal labels
            \node
            at (-1.5708,-0.1) {$-\pi/2$};
                        \node
            at (1.5708,-0.1) {$\pi/2$};

                 \node%[fill=gray!30,inner sep=2pt]
            at (-1.1,0.52) {$H$};
            
                 \node[fill=gray!30,inner sep=2pt]
            at (0.95,0.1) {$S_a$};

        % Table to the right of the plot
        \node[anchor=west,inner sep=0pt] at (2.5,0.24195) {%
            \small
            \renewcommand{\arraystretch}{1}%
            \begin{tabular}{cc}
                \toprule
                \multicolumn{1}{c}{%
                    \shortstack{Partition\\
                                size $n$}}
                &
                \shortstack{Acceptance rate\\
                            $|S_a|/|H|\ (\%)$}
                \\
                \midrule
                10    & 85.80\,\% \\
                50   & 96.55\,\% \\
                100   & 98.25\,\% \\
                500  & 99.64\,\% \\
                1000  & 99.82\,\% \\
                1500  & 99.88\,\% \\
                10000 & 99.99\,\% \\
                \bottomrule
            \end{tabular}%
        };

    \end{tikzpicture}
    \caption{Graph of the angular density with an upper rectangular covering ($n=10$)
    associated with the standard normal distribution
             and acceptance rates in percentages
             for several partition sizes $n$.}
    \label{Normal:rectangular:partition}
\end{figure}

The set $S_a$ as defined in~\eqref{Angular:set:S} is the region below this density.
The plot on the left of
Figure~\ref{Normal:rectangular:partition} shows a rectangular
covering for $S_a$ and the table on the right shows the increase in the acceptance rate
as the partition size $n$ grows. These numbers are the same as those displayed
in the table of Figure~\ref{Normal:circular:partition}.

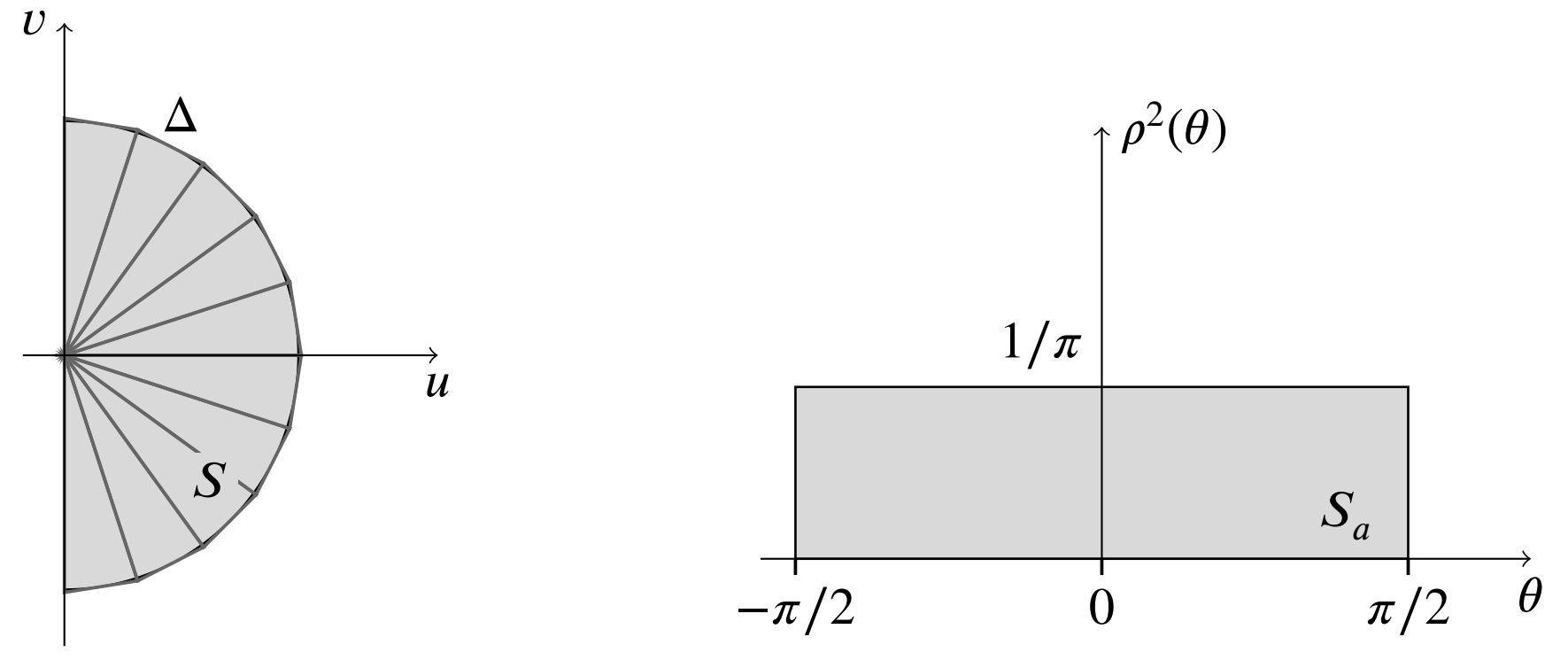
\begin{figure}[H]
    \centering
    \begin{tikzpicture}[xscale=3.1,yscale=3.1]
        \pgfmathsetmacro{\R}{1/sqrt(pi)}
        \pgfmathsetmacro{\h}{1/pi}
        \pgfmathsetmacro{\outerR}{\R/cos(9)}

        % Quotient region: right half-disk
        \filldraw[line width=0.5pt,fill=gray!30]
            plot[
                domain=-90:90,
                samples=181,
                variable=\t
            ]
            ({\R*cos(\t)},{\R*sin(\t)})
            -- cycle;

        % Triangular covering with n=10
        \foreach \i in {0,...,9}{
            \pgfmathsetmacro{\angleA}{-90+18*\i}
            \pgfmathsetmacro{\angleB}{\angleA+18}

            \draw[gray!120,line width=0.7pt]
                (0,0)
                -- ({\outerR*cos(\angleA)},
                    {\outerR*sin(\angleA)})
                -- ({\outerR*cos(\angleB)},
                    {\outerR*sin(\angleB)})
                -- cycle;
        }

        % Coordinate axes
        \draw[->] (0,-0.7) -- (0,0.8)
            node[left] {$v$};
        \draw[->] (-0.1,0) -- (0.9,0)
            node[below] {$u$};
            
        \node%[fill=gray!30,inner sep=2pt]
            at (0.28,0.58) {$\Delta$};

        \node[fill=gray!30,inner sep=2pt]
            at (0.35,-0.3) {$S$};

        \node at (0.4,-0.9)
            {(a) Region $S$.};

        % Second panel
        \begin{scope}[
            shift={(2.5,-0.49)},
            xscale=0.47,
            yscale=1.3
        ]
            % Constant density
            \filldraw[line width=0.5pt,fill=gray!30]
                ({-pi/2},0) rectangle ({pi/2},\h);

            % Coordinate axes
            \draw[->] (-1.75,0) -- (2.2,0)
                node[below] {$\theta$};
            \draw[->] (0,0) -- (0,0.8)
                node[right] {$\rho^2(\theta)$};

            % Horizontal ticks
            \foreach \t/\lab in {
                -1.570796/{-\pi/2},
                 0/{0},
                 1.570796/{\pi/2}
            } {
                \draw[line width=0.7pt]
                    (\t,0) -- (\t,-0.03);
                \node[below]
                    at (\t,-0.02) {$\lab$};
            }

            % Height label
            \node[above left]
                at (0,\h) {$1/\pi$};
                
       \node[fill=gray!30,inner sep=1pt]
            at (1.25,0.08) {$S_a$};

            \node at (0,-0.3)
                {(b) Region $S_a$.};
        \end{scope}
    \end{tikzpicture}
    \caption{Regions $S$ and $S_a$ associated with the
        $\mathrm{Cauchy}(0,1)$ distribution.
        The circular covering $C$ is the region $S$.}
        %covered by ten triangles and its angular density.}
    \label{Cauchy-case}
\end{figure}

As a particularly simple example,
consider the $\mbox{Cauchy}(0,1)$ distribution,
\begin{displaymath}
	h(x) = \frac{1}{\pi(1+x^2)}, \quad x \in \mathds{R}.
\end{displaymath}
In this case $h(x)$ is a normalized probability density.
It follows from~\eqref{Radial:function} that the radial function is constant,
\begin{displaymath}
	\rho(\theta) = \frac{1}{\sqrt{\pi}}, 
	\quad \theta \in \left(-{\pi}/{2}, {\pi}/{2}\right),
\end{displaymath}
and the angular density is $\rho^2(\theta)=1/\pi,\ \theta\in\left(-{\pi}/{2}, {\pi}/{2}\right)$.
Thus, the region $S$ is a semicircle and $S_a$ is a rectangle. These are shown in Figure~\ref{Cauchy-case}.
The left plot shows region $S$ along with a triangular covering
which has acceptance rate
\begin{displaymath}
\frac{|S|}{|\Delta|}=\frac{\pi}{2n}\cot\big(\frac{\pi}{2n}\big).
\end{displaymath}
These acceptance rates in percentages are displayed in the second column of Table~\ref{Cauchy:table} up to three decimal places.
The last value $100\%$ is only a close approximation to that value.
On the other hand, the circular sector covering of $S$  is $S$ itself and the sampling has $100\%$ efficiency.
The right plot of Figure~\ref{Cauchy-case} shows the angular density $\rho^2(\theta)$, which yields
an acceptance rate of $100\%$
for any partition size $n$.

\definecolor{darkpowderblue}{RGB}{0,51,153}

\begin{table}[H]
	\centering
	\begin{tabular}{ccc}
                \toprule
                & & \shortstack{Circular-sector covering of $S$ and}\\
                \multicolumn{1}{c}{%
                    \shortstack{Partition\\
                                size $n$}}
		&  \multicolumn{1}{c}{%
                    \shortstack{Triangular covering of $S$\\
                                Acceptance rate $\ (\%)$}}
		& \multicolumn{1}{c}{%
                    \shortstack{Rectangular covering of $S_a$\\
                                Acceptance rate $\ (\%)$}} \\
		\midrule
		10    & 99.176\% & 100\% \\
		50   & 99.967\% & 100\% \\
		100   & 99.992\% & 100\% \\
		500  & 99.999\% & 100\% \\
		1000  & 99.999\% & 100\% \\
		1500  & 99.999\% & 100\% \\
		10000 & 100\% & 100\% \\
		\bottomrule\smallskip
	\end{tabular}
	\label{Cauchy:table}
	\caption{A comparison of acceptance rates of coverings
	for the $\mbox{Cauchy}(0,1)$ distribution. Acceptance rates are rounded to three decimal places. For the triangular covering,
	all rates are strictly below \(100\%\) for any finite \(n\).}

\end{table}

% -----------------------------------------
\section{Discussion}
% -----------------------------------------

It is given in this work three simple sampling schemes to improve the acceptance rate in the ratio-of-uniforms
method. The idea is to propose coverings composed of $n$ smaller regions with 100\% efficiency in the sampling
and such that the covering can get as close to the intended region as one desires.
It was showed that the acceptance rates of the three sampling schemes satisfy
$A_n(\Delta)\le A_n(C)=A_n(R)$ and that
$\displaystyle\lim_{n\to\infty}A_n(\Delta)=\lim_{n\to\infty}A_n(C)=\lim_{n\to\infty}A_n(R)=1$.\\

For a given partition size $n$, the three  sampling schemes require as initial data
the calculation of the $n$ maxima of the radial function. This is a preparation step and
the increase in the acceptance rates of the sampling schemes is a consequence of this prior work.
For the examples shown, those maxima were calculated numerically using the {\tt R} statistical package,
however, special care should be taken to make sure that the maxima used to construct the coverings
are genuine upper bounds. Once this is done, the algorithms
are as follows: to produce one proposal value for the target density, the algorithms
first require generating a value
of a finite distribution given by the probabilities (or weights) of the $n$ triangles, circular sectors  or rectangles.
Then two independent values of the $\mbox{Unif}(0,1)$ distribution are needed to produce a point $(u,v)$
inside the chosen triangle, circular-sector or rectangle.
For the triangular and circular sector coverings, in case of acceptance, $v/u$ is 
a sample from the target distribution.
For the rectangular covering of the angular density, upon acceptance, $\tan u$ 
is a sample from the target density. In summary, in each sampling scheme three random values are needed to produce a proposal value for the target density.\\

An {\tt R} program was used to calculate the acceptance rates shown in the tables accompanying the examples.
In particular, the {\tt R} command {\tt optimize(f,interval=c(a,b),maximum=TRUE)} was used to calculate the
maximum of a function $f$ over the interval $[a,b]$.
The code is not provided here since it is a straightforward implementation of
the instructions given in the text. While the partitions used in this work were homogeneous, employing problem-specific subdivisions
may further improve acceptance rates.\\

Finally, it is known that the ratio-of-uniforms method can be extended to the case of random
vectors of any finite dimension, including modified versions which introduce additional parameters
to improve the efficiency of the uniform sampling (parameter $r>0$, function $g$), see~\cite{Martino-Luengo-Miguez-2018} and~\cite{Wakefield-Gelfand-Smith-1991}. It remains to be investigated whether the coverings proposed in this work could be extended to those more general versions
of the ratio-of-uniforms method.\\

\section*{Generative AI statement}

During the preparation of this work the authors used the AI tool ChatGPT for grammar corrections, language improvement and assistance with the preparation of figures using TikZ in \LaTeX.
After using this tool/service, the authors reviewed, edited and modified the content as needed
and take full responsibility for the content of the publication.

\section*{Disclosure statement}

No potential conflict of interest was reported by the authors.

\section*{Data availability statement}

No external datasets were used in this study. The numerical results were obtained through simulations implemented in {\tt R}. The corresponding code is available from the authors upon reasonable request.

\section*{Funding}
The author E. D. is grateful to SECIHTI for the graduate scholarship support provided.
%under grant number [Your CVU/Scholarship Number]

% ---------------------------------------

\end{document}